\input amstex
\documentstyle{amsppt}
\magnification=1200
\hsize 15.9truecm
\vsize 564pt
\hfuzz 1pt

\def\C{{\bold C}}
\def\bQ{{\bold Q}}
\def\A{{\Cal A}}
\def\E{{\Cal H}}

\def\O{{\Cal O}}

\def\L{{\Cal L}}
\def\H{{\Lambda}}
\def\M{{\Cal M}}
\def\V{{\Cal V}}
\def\Q{{\Cal Q}}
\def\B{{\Cal B}}
\def\bM{{\overline{\Cal M}}}

\def\cPP{{\overline{\Cal {PP}}}}
\def\PP{\Cal {PP}}
\def\Pol{\operatornamewithlimits{Pol}}
\def\deg{\operatorname{deg}}
\def\Res{\operatornamewithlimits{Res}}
\def\Aut{\operatornamewithlimits{Aut}}
\def\LL{{LL}}
\def\LCM{\operatorname{LCM}}
\def\CP{{\bold CP}}

\topmatter
\title 
On Hurwitz numbers and Hodge integrals
 \endtitle
\author Torsten Ekedahl$^\dag$, Sergei Lando$^\ddag$, 
Michael Shapiro$^\P$, and  Alek Vainshtein$^\S$ 
\endauthor
\headline{\hss{\rm\number\pageno}\hss}
\footline{}
\affil
$^\dag$ Dept. of Math., University of Stockholm, S-10691, Stockholm, 
teke\@matematik.su.se\\ 
$^\ddag$ Math. College, Moscow Independent University, 
Moscow, lando\@lando.mccme.rssi.ru\\ 
$^\P$ Dept. of Math., Royal Institute of Technology, S-10044, 
Stockholm, mshapiro\@math.kth.se\\ 
$^\S$ Dept. of Math. and Dept. of Computer Science, University of Haifa, 
31905 Haifa, alek\@mathcs2.haifa.ac.il
\endaffil
\abstract
In this paper we find an explicit formula for the number of topologically
different ramified coverings $C\to\CP^1$ ($C$ is a compact Riemann surface 
of genus $g$) with only one complicated branching point in terms of Hodge 
integrals over the moduli space of genus $g$ curves with marked points.
\endabstract
\endtopmatter

\document
\subheading{1. Introduction and main results}
For a compact connected genus $g$ complex curve $C$,
let $\lambda\: C\to\CP^1$ be a meromorphic function.
We treat this function as a ramified covering
of the sphere. Two ramified coverings $(C_1;\lambda_1)$,
$(C_2;\lambda_2)$ are called
{\it topologically equivalent\/} if there exists a
homeomorphism $h:C_1\to C_2$ making the following diagram
commutative:
$$
\matrix
C_1& \overset{h}\to\longrightarrow&C_2\\
\lambda_1\searrow && \swarrow\lambda_2\\
&\CP^1
\endmatrix
$$
The critical values of topologically equivalent functions,
i.e., the ramification points of the coverings,
coincide, as do the genera of the covering curves.
In his famous paper~\cite{11}  Hurwitz initiated the
topological classification of such coverings in the case when
exactly one of the ramification points is degenerate, and the
remaining points are nondegenerate.
Below we refer to the degenerate ramification point as ``infinity'',
and its preimages are called ``poles''.
For a given set of orders $k_1,\dots,k_n$ of $n$ distinct poles,
the number of the equivalence classes of
topologically nonequivalent ramified coverings
with these orders of poles and prescribed nondegenerate
ramification points is finite. This number $h_{g;k_1,\dots,k_n}$
(called the {\it Hurwitz number\/}) is independent
of the exact location of the nondegenerate ramification points.
The number of sheets in the covering is $k=k_1+\dots+k_n$,
and the space of all such coverings has dimension $d=k+n+2g-2$
over $\C$.
Hurwitz posed the problem of finding $h_{g;k_1,\dots,k_n}$ explicitly.
There exists a large number of publications on this topic written both by
physicists and mathematicians, see e.g.~\cite{1--3, 8--11, 15, 16}.
Our aim is to express Hurwitz numbers in terms of intersection numbers
for the Chern classes of certain line bundles on the moduli space of
complex curves with $n$ marked points.

Let $\bM_{g;n}$ denote the Deligne--Mumford compactification
of the moduli space $\M_{g;n}$ of genus $g$ curves with $n$ marked points,
and let $\L_i$ be the line bundle on $\bM_{g;n}$
whose fiber above a point $(C;z_1,\dots,z_n)\in\bM_{g;n}$
coincides with the cotangent space to $C$ at $z_i$.
The first Chern class of such a bundle is denoted by $c_1(\L_i)$.

The main result of the present paper is as follows.
Let $\#\Aut(k_1,\dots,k_n)$ denote the number of automorphisms
of the $n$-tuple $(k_1,\dots,k_n)$,
$c(\H_{g;n})$ denote the total Chern class
of the relative dualizing sheaf over $\bM_{g;n}$.

\proclaim{Theorem 1.1}
The Hurwitz number $h_{g;k_1,\dots,k_n}$ equals
$$
\frac{d!}{\#\Aut(k_1,\dots,k_n)}
\prod_{i=1}^n\frac{k_i^{k_i}}{k_i!}
\int_{\bM_{g;n}}\frac{c(\H_{g;n})}{(1-k_1c_1(\L_1))\dots(1-k_nc_1(\L_n))}.
\tag1
$$
\endproclaim
Here we admit as usual that the integral of a class
with the degree different from that of the dimension of the
manifold vanishes.
A sketch of the proof is given in Section~2; a detailed account will appear
elsewhere.

More general expressions
$\int_{\bM_{g;n}} c_{\alpha_1}(\H_{g;n})\cdots c_{\alpha_q}(\H_{g;n})
c_1(\L_1)^{a_1}\cdots c_1(\L_n)^{a_n}$,
including~(1) as a special case, are called {\it Hodge integrals\/};
they are studied extensively in connection with Gromov--Witten invariants,
quantum cohomologies and enumerative geometry
(see e.g.~\cite{4, 7, 12--14, 17}).

Combining Theorem~1.1 with Theorem~6 of~\cite{16},
we recover the following recent calculation by C.~Faber and
R.~Pandharipande~\cite{5}.

\proclaim{Theorem 1.2}
$$
1+\sum_{g=1}^\infty t^{2g}
\int_{\bM_{g;1}}\frac{k^g+k^{g-1}c_1(\H_{g;n})+\dots+c_g(\H_{g;n})}
{1-c_1(\L_1)}=\left(\frac{t/2}{\sin(t/2)} \right)^{k+1}.
$$
\endproclaim

The main ingredient of the proof of Theorem~1.1 is the following
nonlinear
version of the product formula for the total Segre classes of vector bundles.
A similar homogeneous situation is considered in Chapter 5 of~\cite{6}.

Let $\pi\: \A\to X$  be a fiber bundle over a complex $m$-dimensional
variety $X$ whose fibers are quasihomogeneous  spaces
$\C^{r+1}_{w_0,\dots,w_r}$ with weights $(w_0,\dots,w_r)$. Assume for 
simplicity that all $w_i$ are positive integers.
We call $\A$ a {\it bundle of quasihomogeneous spaces\/} if the transition 
functions (possibly not linear) preserve the weights.
Clearly, any vector bundle is a bundle of quasihomogeneous spaces
with weights $(1,\dots,1)$. Any  subbundle $\Q\subset \A$ whose fibers
are quasihomogeneous  cones (with the same weights $(w_0,\dots,w_r)$)
in $\C^{r+1}_{w_0,\dots,w_r}$ is called a
{\it bundle of quasihomogeneous cones}. Note that $\C^*$ acts
fiberwise on both $\A$ and $\Q$ by multiplication of the corresponding
coordinates by $t^{w_i}$. The quotient of $\A$ by the
$\C^*$-action is fibered over $X$, and the fibers are isomorphic to the
weighted projective spaces $P_{w_0,\dots,w_r}$. We denote this quotient
by $P\A$ and say that $p\: P\A\to X$ is a {\it bundle of weighted projective 
spaces}.
The quotient of $\Q$ by the $\C^*$-action is denoted by $P\Q$.
Note that the fibers of both $P\A$ and $P\Q$ can be singular.

Put $N=\LCM(w_0,\dots,w_r)$.
Denote by $\Q\oplus 1$ the fiberwise direct product of $\Q$ with the 
trivial one-dimensional vector bundle. Clearly,  
$\Q\oplus 1$ is also a bundle of quasihomogeneous cones with respect
to the natural embedding $\Q\oplus 1\subset \A\oplus 1$. 
The fibers of $\A\oplus 1$ are just
$\C^{r+2}_{1,w_0,\dots,w_r}$ (the weight of the new coordinate equals 1).
By $q$ we denote the natural projection $q\: \Q\oplus 1\to X$.
There is a natural sheaf $\O(N)$ of quasihomogeneous functions of
quasihomogeneous degree $N$ on $P(\A\oplus 1)$ 
(which restricts also naturally to $P(\A\oplus 1)$).

Following~\cite{6} we define the total Segre class
$s(\Q)\in H^*(X,{\bQ})$ by the formula
$$
s(\Q)=\sum_{i\ge 0} 
q_*\left(\frac{c_1(\O(N))^{i}\cap [P(\Q\oplus 1)]}{N^{i}}\right).
$$

Assume that we are given a quasihomogeneous morphism
$\psi\: \A\to \V$ from a bundle of quasihomogeneous spaces $\A$ into
a vector bundle $\V$  which is fiberwise of quasihomogeneous degree 1.
Consider a {\it nonlinear exact triple}
$$
0\rightarrow\Q_\psi\overset{i}\to\hookrightarrow\A
\overset{\psi}\to\rightarrow\V\to 0,
$$
where the bundle of quasihomogeneous cones
$\Q_\psi$ is the inverse image of the zero section in $\V$.

\proclaim{Theorem 1.3}
Under the above assumptions one has
$$
s(\A)=s(\Q_\psi)\cdot s(\V).
$$
\endproclaim

The authors are sincerely grateful to V.I.Arnold who introduced the last
three authors to this subject and kept us informed about his own research.
We are indebted to B.~Dubrovin, L.~Ernstr\"om, C.~Faber, Yu.~Manin, 
S.~Natanzon, R.~Pandharipande, M.~Rosellen and, especially, B.~Shapiro, 
for many fruitful discussions.

\subheading{2. A sketch of the proof of Theorem~1.1}
The space of all (holomorphic equivalence
classes of) generic meromorphic functions on genus $g$ curves
with prescribed orders $(k_1,\dots,k_n)$
of poles will be called the {\it Hurwitz space\/}
and denoted by $H_{g;k_1,\dots,k_n}$.
It carries a natural complex structure.
Observe that $H_{g;k_1,\dots,k_n}$ is naturally fibered over $\M_{g;n}$:
the projection sends poles to marked points.
There is also a natural mapping
$\LL\: H_{g;k_1,\dots,k_n}\to\Pol_{d}$ called the
{\it Lyashko--Looijenga map}, where $\Pol_{d}$ is the space of monic
polynomials in one variable of degree $d=k+n+2g-2$; $\LL$ associates to each
meromorphic function the unordered set of its critical values or,
equivalently, the polynomial vanishing exactly on this set.

There exists a simple relation between the degree of the Lyashko--Looijenga map
and the corresponding Hurwitz number (cf.~\cite{8}).

\proclaim{Proposition 2.1}
$$
\frac{\deg \LL}{h_{g;k_1,\dots,k_n}}= \#\Aut(k_1,\dots,k_n)\prod_{i=1}^nk_i.
$$
\endproclaim

The following statement plays a key role in our construction.

\proclaim{Proposition 2.2}
There exists a compactification $\E_{g;k_1,\dots,k_n}$ of
$H_{g;k_1,\cdots,k_n}$ which is a bundle over $\bM_{g;n}${\rm ;} moreover, the
Lyashko--Looijenga map extends to $\E_{g;k_1,\dots,k_n}$.
\endproclaim

It is now easy to express $\deg \LL$ in terms of intersection numbers.
Indeed, let
$x^{d}+a_1x^{d-1}+\dots+a_{d}\in \Pol_{d}$; we
denote by $D_i$ the divisor in $\E_{g;k_1,\dots,k_n}$ given by the
quasihomogeneous equation $\{a_i=\text{const}\cdot a_1^i\}$. Then
$\deg \LL=w(k_1,\dots,k_n)\deg ([D_2]\cap\cdots\cap[D_{d}])$,
where $[D_i]\in A_*P\E_{g;k_1,\dots,k_n}$ and $w(k_1,\dots,k_n)$ is defined by
 the quasihomogeneous weights in the image and the preimage. (As usual, given
an algebraic variety $X$, we write $A_* X$ for the rational
equivalence group with coefficients in $\bold Q$; see~\cite{6} for details.)

Put $\Delta=c_1(\Q(1))\in A_*P\E_{g;k_1,\dots,k_n}$.
One can easily show that $[D_i]=i\cdot\Delta$ in the Chow ring. Therefore,
$\deg \LL =d!w(k_1,\dots,k_n)\int_{P\E_{g;k_1,\dots,k_n}}
\Delta^{d}$. By the definition of the top Segre class we have
$\int_{P\E_{g;k_1,\dots,k_n}} \Delta^{d}=
\int_{\bM_{g;n}}s_{top}(\E_{g;k_1,\dots,k_n})=
\int_{\bM_{g;n}}s(\E_{g;k_1,\dots,k_n})$. Thus
$\deg \LL =d!w(k_1,\dots,k_n)\int_{\bM_{g;n}}s(\E_{g;k_1,\dots,k_n})$.

The quasihomogeneity factor $w(k_1,\dots,k_n)$ is given by the following
statement (cf.~\cite{1, 8}).

\proclaim{Proposition 2.3}
$$
w(k_1,\dots,k_n)=\prod_{i=1}^n\frac {k_i^{k_i}}{(k_i-1)!}.
$$
\endproclaim

Denote by $\PP_{g;k_1,\dots,k_n}\to{\M_{g;n}}$ the fiber bundle of
principal parts  of germs of
Laurent polynomials of degrees $k_1,\dots,k_n$ at the marked points
$P_1,\dots,P_n$ on genus $g$ surfaces.

\proclaim{Proposition 2.4}
There exists a compactification $\cPP_{g;k_1,\dots,k_n}$ of
$\PP_{g;k_1,\dots,k_n}$ which is a bundle of quasihomogeneous spaces
over $\bM_{g;n}$; moreover,
$\E_{g;k_1,\dots,k_n}$ is a subbundle of quasihomogeneous cones
in $\PP_{g;k_1,\dots,k_n}$.
\endproclaim

We now define a mapping $f\:\cPP_{g;k_1,\dots,k_n}\to\H_{g;n}$.
Let $p=(p_1,\dots,p_n)\in\cPP_{g;k_1,\dots,k_n}$ be a collection of
principal parts of Laurent germs at $P_1,\dots,P_n$ and $\omega$ be a
holomorphic 1-form.
We put $\langle f(p),\omega\rangle=\sum_{i=1}^n \Res_{P_i} (p_i\omega)$.
Plainly, $\E_{g;k_1,\dots,k_n}$ is the inverse image of the zero section.
Hence, it is included in the nonlinear exact triple
$$
0\rightarrow\E_{g;k_1,\dots,k_n}\overset{i}\to\hookrightarrow
\cPP_{g;k_1,\dots,k_n}\overset{f}\to\rightarrow\H_{g;n}\to 0.
$$
Thus, applying Theorem~1.3 one obtains
$s(\E_{g;k_1,\dots,k_n})={s(\cPP_{g;k_1,\dots,k_n})}/{s(\H_{g;n})}$.
Taking into account that $s(\H_{g;n})=1/{c(\H_{g;n})}$, where $c$ stands
for the total Chern class, we rewrite the above formula as
$s(\E_{g;k_1,\dots,k_n})={s(\cPP_{g;k_1,\dots,k_n})}\cdot{c(\H_{g;n})}$.

On the other hand, projecting a collection of principal parts onto the 
collection 
of its leading terms one gets another nonlinear exact triple:
$$
0\to\B\to\cPP_{g;k_1,\dots,k_n}\to\bigoplus_i(\L_i^*)^{\otimes k_i}\to 0,
$$
where $\B$ is a trivial bundle of quasihomogeneous spaces of nonleading terms
and $\L_i$ are the linear bundles mentioned in the introduction.
Applying Theorem~1.3 one more time we complete the proof of
Theorem~1.1.

\enddocument